\documentstyle[11pt, amsfonts]{amsart}
\begin{document}
\def\R{{\mathbb R}}
\def\Z{{\mathbb Z}}
\newcommand{\trace}{\rm trace}
\newcommand{\Ex}{{\mathbb{E}}}
\newcommand{\Prob}{{\mathbb{P}}}
\newcommand{\E}{{\cal E}}
\newcommand{\F}{{\cal F}}
\newtheorem{df}{Definition}
\newtheorem{theorem}{Theorem}
\newtheorem{lemma}{Lemma}
\newtheorem{pr}{Proposition}
\newtheorem{co}{Corollary}
\def\n{\nu}
\def\sign{\mbox{ sign}}
\def\a{\alpha}
\def\N{{\mathbb N}}
\def\A{{\cal A}}
\def\L{{\cal L}}
\def\X{{\cal X}}
\def\C{{\cal C}}
\def\F{{\cal F}}
\def\c{\bar{c}}
\def\v{\nu}
\def\d{\partial}
\def\diam{\mbox{\rm diam}}
\def\vol{\mbox{\rm Vol}}
\def\b{\beta}
\def\t{\theta}
\def\l{\lambda}
\def\e{\varepsilon}
\def\colon{{:}\;}
\def\pf{\noindent {\bf Proof :  \  }}
\def\endpf{ \begin{flushright}
$ \Box $ \\
\end{flushright}}
\title[The Busemann-Petty problem for arbitrary measures.]
{The Busemann-Petty problem for arbitrary measures.}
\begin{abstract}
The aim of this paper is to study properties of sections of convex bodies with 
respect to different types of measures. We present a formula 
connecting the Minkowski functional of a convex symmetric body $K$ with the
 measure of its sections. We apply this formula to study properties of
general measures most of  which were known before only in the case of the 
standard Lebesgue measure. We solve an analog of the Busemann-Petty problem 
for the case of general measures. In addition, we show that there are measures,
 for 
which the answer to the generalized Busemann-Petty problem is affirmative 
in all dimensions. Finally, we apply the latter fact to prove a number of 
different inequalities concerning the volume of sections of convex symmetric
 bodies in $\R^n$ and solve a version of  generalized Busemann-Petty problem
for sections by $k$-dimensional subspaces.
\end{abstract}

\author{A. Zvavitch}
\address{
Artem Zvavitch\\
Department of Mathematical Sciences\\
Kent State University\\
Kent, OH 44242, USA}

\email{zvavitch@@math.kent.edu}

\maketitle
\section{Introduction}
\noindent Consider a non-negative, even function $f_n(x)$, which is locally 
integrable on $\R^n$. Let $\mu_n$ be the measure on $\R^n$ with  density
$f_n$. 

For $\xi\in S^{n-1}$, let $\xi^\bot$ be the central hyperplane orthogonal to 
$\xi$. Define a measure $\mu_{n-1}$ on $\xi^\bot$, for each $\xi\in S^{n-1}$, 
so that for every bounded Borel set $B\subset\xi^{\bot}$,
$$
\mu_{n-1}(B)=\int\limits_B f_{n-1}(x)dx,
$$
where $f_{n-1}$ is an even function on $\R^n$, which is locally integrable on
each $\xi^\bot$.

In this paper we study the following problem

\noindent{\bf The Busemann-Petty problem for general measures (BPGM): }

\noindent {\it Fix $n\ge 2$. Given two
convex origin-symmetric bodies $K$ and $L$ in $\R^n$ such that
$$
\mu_{n-1}(K\cap \xi^\bot)\le\mu_{n-1}(L\cap \xi^\bot) 
$$
for every $\xi \in S^{n-1}$, does it follow that
$$
\mu_n(K)\le \mu_n(L)?
$$}
Clearly, the BPGM problem is a triviality for $n=2$ and $f_{n-1}>0$,
 and the answer is ``yes'',  moreover $K \subseteq L$.
Also note that this problem is a generalization of the  Busemann-Petty  problem,  
posed in 1956 (see \cite{BP}) and   asking the same question for
 Lebesgue measures: $\mu_n(K)=\vol_n(K)$ and 
$\mu_{n-1}(K\cap\xi^\bot)=\vol_{n-1}(K\cap\xi^\bot)$ i.e. 
$f_n(x)=f_{n-1}(x)=1$.

Minkowski's theorem (see \cite{Ga3}) shows that an origin-symmetric
 star-shaped body is uniquely determined by the volume of its 
hyperplane sections (the same is true for the case of general 
symmetric measure,
 see Corollary \ref{co:un} below). In view of this fact it is 
quite surprising that the answer to the original Busemann Petty problem is 
negative for 
$n\ge 5$. Indeed, it is affirmative if $n\le 4$ and negative if 
$n\ge 5$. The solution appeared as the result   
of a sequence of papers: \cite{LR} $n\ge 12$, \cite{Ba} $n\ge 10$, \cite{Gi} and \cite{Bo2} $n\ge 7$,  \cite{Pa} and \cite{Ga1} $n\ge 5$,  
\cite{Ga2} $n=3$, \cite{Zh2} and \cite{GKS} $n=4$   
(we refer to \cite{Zh2}, \cite{GKS} and \cite{K8} for more historical details).

 It was  
shown in \cite{Z}, that the answer to BPGM  in the case of the 
standard Gaussian measure ($f_n(x)=f_{n-1}=e^{-|x|^2/2}$) is the same:
 affirmative if $n\le 4$ and negative if $n \ge 5$.

Answers to the original and Gaussian Busemann-Petty problems suggest that  
the answer to the BPGM could be independent from the choice of measures and depend only on the dimension $n$.

In  Corollary \ref{co:me} below we confirm this conjecture by proving the
following:
\vskip .3cm
\noindent{\it Let  $f_n(x)=f_{n-1}(x)$ be equal even nonnegative continuous functions, 
then the answer to the  BPGM problem
 is affirmative  if $n\le 4$ and negative if $n\ge 5.$}
\vskip .3cm
Actually, the above fact is a corollary of a pair of more general theorems. 
Those theorems use the Fourier transform in the sense of distributions to 
characterize those functions $f_n(x)$ and $f_{n-1}(x)$ for 
which  the BPGM problem has affirmative (or negative) answer in a given 
dimension:
\begin{theorem}\label{th:maina} {\bf (BPGM: affirmative case)}
Let $f_n$ and $f_{n-1}$ be even continuous nonnegative functions such that
\begin{equation}\label{increas}
t\frac{f_n(tx)}{f_{n-1}(tx)}
\end{equation}
is an increasing function of $t$ for any fixed $x\in S^{n-1}$.
Consider a symmetric star-shaped body $K$  in $\R^n$ such that
\begin{equation}\label{eq:pos}
\|x\|_K^{-1}\frac{f_n\left(\frac{x}{\|x\|_K}\right)}
{f_{n-1}\left(\frac{x}{\|x\|_K}\right)}
\end{equation}
is a positive definite distribution on $\R^n$.
Then for any symmetric  star-shaped body $L$ in $\R^n$ satisfying
$$
\mu_{n-1}(K\cap \xi^\bot)\le \mu_{n-1}(L\cap \xi^\bot),\,\,\,\forall \xi\in 
S^{n-1},
$$
we have
$$
\mu_n(K)\le \mu_n(L).
$$
\end{theorem}
\begin{theorem}\label{th:mainn}{\bf (BPGM: negative case)} Let $f_n$ and 
$f_{n-1}$ be even continuous nonnegative functions such that
$$
t\frac{f_n(tx)}{f_{n-1}(tx)}
$$
is an increasing function of $t$ for any fixed $x\in S^{n-1}$. Also
 assume that $f_{n-1}(x)\in C^\infty(\R^n\setminus \{0\})$ and is it strictly 
positive on  $\R^n\setminus \{0\}$.
If $L$ is an infinitely smooth, origin symmetric, convex body in $\R^n$ 
with positive curvature, and the function
\begin{equation}\label{con2}
\|x\|_L^{-1}\frac{f_n\left(\frac{x}{\|x\|_L}\right)}
{f_{n-1}\left(\frac{x}{\|x\|_L}\right)}
\end{equation}
 is in $C^\infty(\R^n\setminus \{0\})$ and 
does not represent a positive definite distribution,
 then there exists a convex symmetric body $D$ in $\R^n$ such that
$$
\mu_{n-1}(D\cap \xi^\bot)\le \mu_{n-1}(L\cap \xi^\bot),\,\,\,\forall \xi\in
 S^{n-1},
$$
but
$$
\mu_n(D) > \mu_n(L).
$$
\end{theorem}
Note that the differentiability assumptions in Theorem \ref{th:mainn} and
the assumption on $f_{n-1}$ to be strictly positive are
not critical for  most of applications and  can be avoided by applying a 
standard approximation argument (see Section \ref{sec:apl}).

Theorems \ref{th:maina} and \ref{th:mainn} are generalizations of a
theorem of Lutwak (see \cite{Lu}) who provided a characterization of symmetric
star-shaped bodies for which the original Busemann-Petty problem has an 
affirmative answer (see \cite{Z} for the case of Gaussian measure). 
Let $K$ and $M$ be symmetric star-shaped bodies in $\R^n.$ We say that $K$ is
the intersection body of $M$ if the radius of $K$ in every direction
is equal to the $(n-1)$-dimensional volume of the central hyperplane 
section of $L$ perpendicular to this direction. A more general class of 
{\it intersection bodies} is defined as the closure in the radial metric
 of the class of intersection bodies of star-shaped  bodies (see \cite{Ga3},
 Chapter 8).

Lutwak (\cite{Lu}, see also \cite{Ga2} and \cite{Zh1}) proved that if $K$ is 
an intersection body then 
the answer to the original Busemann-Petty  problem is affirmative for every 
$L$, and, on the other hand, if $L$ is 
 not an intersection body, then one can perturb 
it to construct a body $D$ giving together with $L$ a counterexample. 

Lutwak's result is related to Theorems  \ref{th:maina} and \ref{th:mainn}
via the  following Fourier analytic characterization of intersection bodies 
found by  Koldobsky  \cite{K3}: 
an origin symmetric star body $K$ in $\R^n$ is an intersection body if
and only if the function $\|\cdot\|_K^{-1}$ represents a 
positive definite distribution on $\R^n.$

We present the proof of Theorems \ref{th:maina} and \ref{th:mainn}
 in Section \ref{proof}.
The proof is based on the Fourier transform 
of distributions, the Spherical Parseval's 
identity introduced by Koldobsky (see Lemma 3 in \cite{K4} or Proposition 
\ref{pro} in Section \ref{proof}) and an elementary functional
inequality (see Lemma \ref{ell}).

Another application of Theorem \ref{th:maina} 
is motivated by a question of what one has to know about the measure 
of central sections of the bodies $K$ and $L$ to make a conclusion about the
relation between the volumes of $K$ and $L$ in every dimension. Results
of such a  type, involving derivatives or the Laplace transform 
of the parallel sections functions were proved in \cite{K4}, \cite{K6},
 \cite{K7}, \cite{K8}, \cite{K9}, \cite{RZ}, \cite{KYY}.

 Note that Theorem \ref{th:maina} allows us to start a  different approach 
to this problem which is based on introducing of two  different measures: 
$\mu_n$ on convex bodies and $\mu_{n-1}$ on  hyperplane sections of 
convex bodies in $\R^n$. This  leads to a number of interesting facts and 
gives  examples of non-trivial  densities $f_n(x)$ and $f_{n-1}(x)$ 
for which the BPGM problem has an affirmative answer in any dimension (see Section
\ref{sec:apl}). Probably the most notable statement is 
(see Corollary \ref{co:p}):
\vskip .3cm
\noindent{\it For any $n\ge 2$ and any symmetric star-shaped bodies 
$K, L\subset \R^n$ such that
$$
\int_{K\cap \xi^\bot}\sum\limits_{i=1}^n|x_i|\, dx \le \int_{L\cap \xi^\bot}
\sum\limits_{i=1}^n|x_i|\, dx 
$$ 
for every $\xi\in S^{n-1}$, we have
$$  
\vol_n(K)\le \vol_n(L).  
$$  
}
\vskip .3cm
One of the advantages of the latter result (and other results of this type) is
that now one can apply methods and inequalities from
asymptotic  convex geometry (see \cite{MP}) to produce new bounds on the volume
 of hyperplane sections of a convex body (see Section \ref{sec:apl}). 

Finally iterating Theorem  \ref{th:maina} with  different densities $f_n$
 and  $f_{n-1}$ we present some new results for sections of
 codimension
greater than $1$ (see Theorem \ref{th:codim} and Corollary \ref{co:krz}).


\section{Measure of sections of star-shaped bodies and the Fourier transform}\label{sec:f}


Our main tool is the Fourier transform of distributions (see \cite{GS},
 \cite{GV} and \cite{K8} for exact definitions and properties). 
We denote by ${\cal S}$ the space of rapidly decreasing 
infinitely differentiable functions (test functions) on $\R^n$ with values in 
${\mathbb C}$. By ${\cal S}'$ we denote the space of distributions over 
$\cal S$. The Fourier transform of a distribution $f$ is defined by 
$ \langle \hat{f},
\hat{\phi}\rangle=(2\pi)^n\langle f,\phi\rangle,$ 
for every test function $\phi$. A distribution $f$ is called {\it even 
homogeneous} of degree $p \in \R$ if 
$$
\langle f(x), \phi(x/t)\rangle=|t|^{n+p}\langle f(x), \phi(x)\rangle,\,\,\, 
\forall \phi\in {\cal S}, \,\, t\in \R\setminus\{0\}.
$$
The Fourier transform of an even homogeneous distribution of degree $p$ is 
an even homogeneous distribution of degree $-n-p$.

A distribution $f$ is called {\it positive definite} if, for every nonnegative 
test function  $\phi\in S$,
$$
\langle f,\phi*\overline{\phi(-x)}\rangle\ge 0.
$$
By Schwartz's generalization of Bochner's theorem, a distribution is positive
 definite if and only if its Fourier transform is a positive distribution 
(in the sense that $\langle \hat{f},\phi\rangle\ge 0$, for every non-negative 
$\phi\in S$). 
Every positive distribution is a tempered measure, i.e. a Borel non-negative, 
locally finite measure $\gamma$ on $\R^n$ such that, for some $\beta>0$,
$$
\int_{\R^n}(1+|x|)^{-\beta}d\gamma(x) <\infty,
$$
where $|x|$ stands for the Euclidean norm (see \cite{GV} p. 147).

The {\it spherical Radon transform} is a bounded linear operator on 
$C(S^{n-1})$ defined by
$$
{\cal R} f(\xi)=\int_{S^{n-1}\cap \xi^\bot} f(x) dx,\,\,\, f\in C(S^{n-1}),\,\, \xi\in S^{n-1}.
$$
Koldobsky (\cite{K1}, Lemma 4) proved that if $g(x)$ is an even homogeneous 
function of degree $-n+1$ on $\R^n\setminus \{0\}$, $n>1$ so that 
$g\big|_{S^{n-1}}\in L_1(S^{n-1})$, then 
\begin{equation}\label{eq:kol}
{\cal R}g(\xi)=\frac{1}{\pi}\hat{g}(\xi), \,\,\,\, \,\,\,\,\, \forall \xi \in S^{n-1}.
\end{equation}
Let $K$ be a body (compact set, with non-empty interior) that is star-shaped 
with respect to the origin in $\R^n$. The Minkowski functional of $K$ is given
 by
$$
\|x\|_K=\min\{\alpha>0: x\in \alpha K\}, \,\,\, x\in \R^n.
$$
\begin{theorem}\label{tm:f} Let $K$ be a  symmetric star-shaped body in $\R^n$, then 
$$
\mu_{n-1}(K\cap \xi^\bot)=\frac{1}{\pi}\left(|x|^{-n+1}\int\limits_{0}^{|x|/\|x\|_K} t^{n-2} f_{n-1}\left(\frac{tx}{|x|}\right) dt\right)^{\wedge}(\xi).
$$
\end{theorem}
\pf 
If $\chi$ is the indicator function of the interval $[-1,1]$ then, passing to the polar coordinates in the hyperplane $\xi^\bot$ we get
$$
\mu_{n-1}(K\cap\xi^\bot)=\int\limits_{(x,\xi)=0}\chi(\|x\|_K) f_{n-1}(x)dx
=\int\limits_{S^{n-1}\cap\xi^\bot}\int\limits_0^{\|\theta\|^{-1}_K} t^{n-2}
 f_{n-1}(t\theta)dt d\theta.
$$
We extend the function under the spherical integral  to a homogeneous  of degree $-n+1$ function on $\R^n$ and apply   (\ref{eq:kol}) to get
\begin{eqnarray}\label{e:gh}
\mu_{n-1}(K\cap\xi^\bot)&=&\int\limits_{S^{n-1}\cap\xi^\bot}|x|^{-n+1}\int\limits_0^{\frac{|x|}{\|x\|_K}} t^{n-2} f_{n-1}\left(\frac{tx}{|x|}\right)dt dx\nonumber\\
&=&
{\cal R }\left( |x|^{-n+1}\!\!\!\int\limits_0^{\frac{|x|}{\|x\|_K}} t^{n-2} 
f_{n-1}\left(\frac{tx}{|x|}\right)dt \right)(\xi)\nonumber\\
&=&\frac{1}{\pi}\left( |x|^{-n+1}\int\limits_0^{|x|/\|x\|_K} t^{n-2} 
f_{n-1}\left(\frac{tx}{|x|}\right)dt \right)^\wedge(\xi).
\end{eqnarray}
\endpf
Theorem \ref{tm:f}  implies that a symmetric star-shaped body is uniquely
 determined by the  measure $\mu_{n-1}$ of its sections:   
\begin{co}\label{co:un}
Assume $f_{n-1}(x)\not=0$ everywhere except for a countable set of points in
$\R^n$. Let $K$ and $L$ be star-shaped origin symmetric bodies in $\R^n$. If  
$$  
\mu_{n-1}(K\cap \xi^\bot)=\mu_{n-1}(L\cap \xi^\bot),\,\,\,\,\,\,
 \forall \xi \in S^{n-1},
$$
then $K=L$.  
\end{co}  
\pf    Note that the function in (\ref{e:gh})  
is  homogeneous of degree $-1$ (with respect to $\xi\in \R^n$).  
 This gives a natural extension of $\mu_{n-1}(K\cap\xi^\bot)$ to
 a homogeneous function of degree $-1$. So from the equality of functions
$\mu_{n-1}(K\cap \xi^\bot)=\mu_{n-1}(L\cap \xi^\bot)$ on $S^{n-1}$ we
get the equality of those functions on $\R^n$:
$$
\left(\!|x|^{-n+1}\!\!\int\limits_0^{\frac{|x|}{\|x\|_K}} \!\!\!t^{n-2} 
f_{n-1}\left(\frac{tx}{|x|}\right)dt \right)^\wedge\!\!\!(\xi)=
\left(\!|x|^{-n+1}\!\!\int\limits_0^{\frac{|x|}{\|x\|_L}}\!\!\! t^{n-2} 
f_{n-1}\left(\frac{tx}{|x|}\right)dt \right)^\wedge\!\!\!(\xi).
$$ 
Applying the inverse Fourier transform to  both sides of the latter equation 
we get:  
$$  
\int\limits_{0}^{|x|/||x||_K} t^{n-2}f_{n-1}\left(\frac{tx}{|x|}\right) dt=  
\int\limits_{0}^{|x|/||x||_L} t^{n-2} f_{n-1}\left(\frac{tx}{|x|}\right)dt, \,\,\,\, \forall x\in\R^n,  
$$  
which, together with monotonicity of the function 
$\int\limits_{0}^{y} t^{n-2}f_{n-1}\left(\frac{tx}{|x|}\right) dt$,
 for $y\in\R^+$ (because $f_{n-1}(x)>0$ everywhere except for a countable
 set of points $x\in \R^n$),  gives  $||x||_K=||x||_L$.  
\endpf  


\section{Proofs of Theorems \ref{th:maina} and \ref{th:mainn}}\label{proof}

We would like to start with  the following elementary  inequality:
\begin{lemma}\label{ell}{\bf (Elementary inequality)} 
\begin{align}\label{e:ell}
\int\limits_0^a t^{n-1} \alpha(t)dt - a\frac{\alpha(a)}{\beta(a)}& \int\limits_{0}^a t^{n-2} \beta(t)dt \nonumber \\
& \le \int\limits_0^b t^{n-1} \alpha(t)dt -a\frac{\alpha(a)}{\beta(a)}\int\limits_{0}^b t^{n-2} \beta(t) 
 dt.
\end{align}
for all $a,b > 0$ and $\alpha(t), \beta(t)$  being nonnegative 
 functions on $(0, \max\{a,b\}]$, 
such that all integrals in (\ref{e:ell}) are defined and  $t\frac{\alpha(t)}
{\beta(t)}$ is increasing on $(0, \max\{a,b\}]$.
\end{lemma}
\pf The  inequality (\ref{e:ell}) is equivalent to
$$
a\frac{\alpha(a)}{\beta(a)}  \int\limits_{a}^b t^{n-2} \beta(t)dt
\le \int\limits_a^b t^{n-1} \alpha(t)dt.
$$
But
$$
a\frac{\alpha(a)}{\beta(a)}  \int\limits_{a}^b t^{n-2} \beta(t)dt=
\int\limits_{a}^b  t^{n-1} \alpha(t) \left(a\frac{\alpha(a)}{\beta(a)}\right)
\left(t\frac{\alpha(t)}
{\beta(t)}\right)^{-1}dt
\le \int\limits_a^b t^{n-1} \alpha(t)dt.
$$
Note that the latter inequality does not require $a\le b$.
\endpf
Before proving  Theorem \ref{th:maina} (the affirmative case of BPGM) we 
need to state a version of Parseval's identity on the sphere and to remind 
a few facts concerning  positive definite homogeneous distributions.

Suppose that $f(x)$ is a continuous on $\R^n\setminus\{0\}$ function, which  
is a positive definite homogeneous distribution of degree $-1$. Then the 
Fourier transform of $f(x)$ is a tempered measure $\gamma$ on $\R^n$ (see 
Section  \ref{sec:f}) which is  a homogeneous distribution of degree 
$-n+1$. Writing this measure in the spherical coordinates (see \cite{K5},
 Lemma 1) we can find a measure $\gamma_0$ on $S^{n-1}$ so that for every 
even test function $\phi$
$$
\langle \hat{f}, \phi\rangle=\langle \gamma, \phi\rangle=
\int_{S^{n-1}}d\gamma_0(\theta)\int_0^\infty \phi(r\theta)dr.
$$

\begin{pr}\label{pro}{\bf (Koldobsky, \cite{K4})}
Let $f$ and $g$ be two functions on $\R^n$,   continuous on $S^{n-1}$ and homogeneous of degrees $-1$ and 
$-n+1$,  respectively. Suppose that  $f$ represents a positive definite
distribution and $\gamma_0$ is the measure on $S^{n-1}$ defined above. Then
$$
\int\limits_{S^{n-1}}\hat{g} (\theta)d\gamma_0(\theta)=(2\pi)^n
\int\limits_{S^{n-1}}f(\theta)g (\theta)d\theta.
$$
\end{pr}
\noindent{\bf Remark:} It is crucial that the sum of degrees  of homogeneity 
of the functions $f$ and $g$   is equal to
$-n$. This is one of the reasons for the choice of degrees of homogeneity in
 the conditions (\ref{eq:pos}) and (\ref{con2}) from Theorems \ref{th:maina},
 and \ref{th:mainn} and in the
formula from Theorem \ref{tm:f}.

\noindent{\bf Proof of Theorem \ref{th:maina}:} Consider symmetric star-shaped bodies $K$ and $L$ in $\R^n$, such that
\begin{equation}\label{e:bn}
\mu_{n-1}(K\cap \xi^\bot)\le \mu_{n-1}(L\cap \xi^\bot),\,\,\,\forall \xi\in S^{n-1}.
\end{equation}
 We apply Theorem \ref{tm:f} to get an analytic form of 
 (\ref{e:bn}):
$$
\left(\!|x|^{-n+1}\!\!\int\limits_0^{\frac{|x|}{\|x\|_K}} \!\!\!t^{n-2} 
f_{n-1}\left(\frac{tx}{|x|}\right)dt \right)^\wedge\!\!\!(\xi)\le
\left(\!|x|^{-n+1}\!\!\int\limits_0^{\frac{|x|}{\|x\|_L}}\!\!\! t^{n-2} 
f_{n-1}\left(\frac{tx}{|x|}\right)dt \right)^\wedge\!\!\!(\xi).
$$
Next we integrate the latter inequality over $S^{n-1}$ with respect to the
measure $\gamma_0$ corresponding to a positive definite homogeneous of degree
$-1$ distribution (\ref{eq:pos}):
\begin{align}
\int\limits_{S^{n-1}}&\left(|x|^{-n+1}\int\limits_0^{\frac{|x|}{\|x\|_K}} t^{n-2} 
f_{n-1}\left(\frac{tx}{|x|}\right)dt \right)^\wedge(\xi) \,d \gamma_0(\xi) \nonumber\\
&\le\int\limits_{S^{n-1}}\left(|x|^{-n+1}\int\limits_0^{\frac{|x|}{\|x\|_L}} t^{n-2} 
f_{n-1}\left(\frac{tx}{|x|}\right)dt \right)^\wedge(\xi)\,d\gamma_0(\xi).
\end{align}
Applying the spherical  Parseval identity (Proposition 
\ref{pro}) we get:
\begin{align}
\int\limits_{S^{n-1}}\|x\|_K^{-1}&\frac{f_n(\frac{x}{\|x\|_K})}
{f_{n-1}(\frac{x}{\|x\|_K})}\int\limits_{0}^{\|x\|^{-1}_K} t^{n-2}
 f_{n-1}(tx)dt d x \nonumber\\
\le&\!\!\int\limits_{S^{n-1}}\|x\|_K^{-1}\frac{f_n(\frac{x}{\|x\|_K})}
{f_{n-1}(\frac{x}{\|x\|_K})}\int\limits_{0}^{\|x\|^{-1}_L} t^{n-2} f_{n-1}(tx) 
\label{1}
 dt d x.
\end{align}
Now we  apply Lemma \ref{ell}, with $a= \|x\|_K^{-1}$,  $b= \|x\|_L^{-1}$,
 $\alpha(t)=f_n(tx)$ and $\beta(t)=f_{n-1}(tx)$ (note that from condition
(\ref{increas}) it follows that $t\alpha(t)/\beta(t)$ is increasing)
 to get
\begin{align}
\int\limits_0^{\|x\|_K^{-1}}&t^{n-1} f_n(tx)dt - \|x\|_K^{-1}
\frac{f_n(\frac{x}{\|x\|_K})}{f_{n-1}(\frac{x}{\|x\|_K})}
 \int\limits_{0}^{\|x\|_K^{-1}} t^{n-2} f_{n-1}(tx)dt\nonumber\\
\le& \int\limits_0^{\|x\|_L^{-1}}\!\! t^{n-1} f_n(tx)dt -
\|x\|_K^{-1}\frac{f_n(\frac{x}{\|x\|_K})}{f_{n-1}(\frac{x}{\|x\|_K})}
 \int\limits_{0}^{\|x\|_L^{-1}}\!\! t^{n-2} f_{n-1}(tx) dt, 
\, \forall x\in S^{n-1}.\nonumber
\end{align}
Integrating over $S^{n-1}$ we get
\begin{align}\label{2}
&\int\limits_{S^{n-1}}\!\int\limits_0^{\|x\|_K^{-1}}\!\!t^{n-1} f_n(tx)dtdx  
- \!\!\!\int\limits_{S^{n-1}}\!\|x\|_K^{-1} \frac{f_n(\frac{x}{\|x\|_K})}
{f_{n-1}(\frac{x}{\|x\|_K})}\int\limits_{0}^{\|x\|_K^{-1}}\!\!t^{n-2}
 f_{n-1}(tx)dtdx\\
&\le \int\limits_{S^{n-1}}\int\limits_0^{\|x\|_L^{-1}}\!\!t^{n-1} f_n(tx)dt dx-
\int\limits_{S^{n-1}} \|x\|_K^{-1} \frac{f_n(\frac{x}{\|x\|_K})}
{f_{n-1}(\frac{x}{\|x\|_K})}\int\limits_{0}^{\|x\|_L^{-1}}\!\!t^{n-2} 
f_{n-1}(tx)  dtdx. \nonumber
\end{align}
Adding equations (\ref{1}) and (\ref{2}) we get
$$
\int\limits_{S^{n-1}}\int\limits_0^{\|x\|^{-1}_K} t^{n-1} f_n(tx)dt dx
\le \int\limits_{S^{n-1}}\int\limits_0^{\|x\|^{-1}_L} t^{n-1} f_n(tx)dt dx,
$$
which is exactly $\mu_n(K)\le \mu_n(L)$.
\endpf
The next proposition is to show that convexity is preserved under small 
perturbation, which is needed in the proof of Theorem  \ref{th:mainn}.
\begin{pr}\label{pr:c}
Consider an infinitely smooth  origin symmetric convex body $L$ with positive 
curvature and even functions $f_{n-1},g \in C^2(R^n\setminus\{0\})$,
such that $f_{n-1}$ is strictly positive on  $\R^n\setminus\{0\}$.
For $\e>0$ define a star-shaped body $D$ by an equation for its radial
 function 
$\|x\|_D^{-1}$:
$$
\int\limits_{0}^{\|x\|_D^{-1}} t^{n-2} 
f_{n-1}
(tx)dt=\int\limits_{0}^{\|x\|_L^{-1}} t^{n-2} f_{n-1}
(tx)dt- \e g(x),\,\,\, \forall x\in S^{n-1}.
$$
Then if $\e$ is small enough the body $D$ is convex.
\end{pr}
\pf  For small enough $\e$,  define a function
 $\alpha_\e(x)$ on $S^{n-1}$ such that
\begin{equation}\label{eq:a}
\int\limits_{0}^{||x||_L^{-1}}\!\!t^{n-2} f_{n-1}
(tx)dt- \e g(x)=\!\!\!
\int\limits_{0}^{||x||_L^{-1}-\alpha_\e(x)}\!\! t^{n-2} f_{n-1}
(tx)dt, \,\,\forall x\in S^{n-1}\!\!.
\end{equation}
Using monotonicity of $\int\limits_{0}^{y} t^{n-2} f_{n-1}(tx)dt$, for $y\in \R^+$ ( $f_{n-1}(tx)>0$, for $tx \in \R^n\setminus\{0\}$), we get
\begin{equation}\label{eq:d}  
||x||_D^{-1}=||x||^{-1}_L-\alpha_\e(x),\,\,\, \forall x\in S^{n-1}. 
\end{equation}
Moreover, using that $f_{n-1}(x)$, $x\in S^{n-1}$,  and its partial derivatives of order
 one are bounded for finite values of $x$ we get, from (\ref{eq:a}), that
 $\alpha_\e(x)$ and its first and second derivatives
converge uniformly to $0$ (for  $x\in S^{n-1}$, as $\e \to 0$).   
Using that $L$ is convex with positive   
curvature, one can choose a small enough $\e$ so that   
the body $D$ is convex (with positive curvature).
\endpf
Next we would like to remind a functional version (with an additional 
differentiability assumption) of Proposition \ref{pro}.
\begin{pr}\label{procont}{\bf (Koldobsky, \cite{K4})} 
Let $f$ and $g$ be two homogeneous $C^\infty(\R^n\setminus\{0\})$ functions of degrees $-1$ and 
$-n+1$,  respectively, then
$$
\int\limits_{S^{n-1}}\hat{f} (\theta)\hat{g} (\theta)d\theta=(2\pi)^n
\int\limits_{S^{n-1}}f(\theta)g (\theta)d\theta.
$$
\end{pr}
We will also need the following fact, which follows from Theorem 1 in 
\cite{GKS}:
if $f$ is positive, symmetric, homogeneous function of degree $-1$, such that $f(x)\in 
C^{\infty}(\R^n\setminus\{0\})$, then $\hat{f}(x)$ is also an infinitely smooth
 function on the sphere $S^{n-1}$.

\noindent{\bf Proof of Theorem \ref{th:mainn}:} First we will use that function (\ref{con2}) is in $C^{\infty}(\R^n\setminus\{0\})$, to claim that
$$
\left(\|x\|_L^{-1}\frac{f_n\left(\frac{x}{\|x\|_L}\right)}
{f_{n-1}\left(\frac{x}{\|x\|_L}\right)}\right)^\wedge
$$
is a continuous function on $S^{n-1}$.

This function does not represent a positive definite distribution so it must be negative on some open symmetric subset $\Omega$ of $S^{n-1}$. Consider  non-negative 
 even function supported $h\in C^\infty(S^{n-1})$ in $\Omega$. Extend $h$ to a homogeneous function 
$h(\theta)r^{-1}$ of degree $-1$. Then the Fourier transform of $h$ is a
homogeneous function of degree $-n+1$: $\widehat{h(\theta)r^{-1}}=
g(\theta)r^{-n+1}$.

For $\e>0$, we  define a body $D$ by
\begin{align}
|x|^{-n+1}\int\limits_{0}^{\frac{|x|}{\|x\|_D}} t^{n-2} &f_{n-1}
\left(\frac{tx}{|x|}\right)dt\nonumber\\
&=|x|^{-n+1}\int\limits_{0}^{\frac{|x|}{\|x\|_L}} t^{n-2} f_{n-1}
\left(\frac{tx}{|x|}\right)dt- \e g\left(\frac{x}{|x|}\right)|x|^{-n+1}.
\nonumber
\end{align}
By Proposition \ref{pr:c} one can choose a small enough $\e$ so that the body $D$ is convex.

Since $h\ge 0$, we have
$$
\mu_{n-1}(D\cap \xi^\bot)=\frac{1}{\pi}\left(|x|^{-n+1}\int\limits_{0}
^{|x|/||x||_D} t^{n-2} f_{n-1}\left(\frac{tx}{|x|}\right)dt\right)^\wedge
(\xi)
$$
$$
=\frac{1}{\pi}\left(|x|^{-n+1}\int\limits_{0}^{|x|/||x||_L} t^{n-2}
f_{n-1}\left(\frac{tx}{|x|}\right)dt\right)^\wedge (\xi)
-\frac{(2\pi)^n\e h(\xi)}{\pi}
$$
$$
\le \mu_{n-1}(L\cap \xi^\bot).
$$
On the other hand, the function $h$ is positive only where
$$
\left(\|x\|_L^{-1}\frac{f_n\left(\frac{x}{\|x\|_L}\right)}
{f_{n-1}\left(\frac{x}{\|x\|_L}\right)}\right)^\wedge(\xi)
$$
is negative so
$$
\left(\|x\|_L^{-1}\frac{f_n\left(\frac{x}{\|x\|_L}\right)}
{f_{n-1}\left(\frac{x}{\|x\|_L}\right)}\right)^\wedge(\xi)
\left(|x|^{-n+1}\int\limits_{0}^{|x|/||x||_D} t^{n-2} f_{n-1}\left(\frac{tx}{|x|}\right)dt\right)^\wedge(\xi)
$$
$$
=\left(\|x\|_L^{-1}\frac{f_n\left(\frac{x}{\|x\|_L}\right)}
{f_{n-1}\left(\frac{x}{\|x\|_L}\right)}\right)^\wedge(\xi)
\left(|x|^{-n+1}\int\limits_{0}^{|x|/||x||_L} t^{n-2} f_{n-1}\left(\frac{tx}{|x|}\right) dt\right)^\wedge(\xi)
$$
$$- (2\pi)^n\left(\|x\|_L^{-1}\frac{f_n\left(\frac{x}{\|x\|_L}\right)}
{f_{n-1}\left(\frac{x}{\|x\|_L}\right)}\right)^\wedge(\xi)
 \e h(\xi)
$$
$$
\ge  \left(\|x\|_L^{-1}\frac{f_n\left(\frac{x}{\|x\|_L}\right)}
{f_{n-1}\left(\frac{x}{\|x\|_L}\right)}\right)^\wedge(\xi)
\left(|x|^{-n+1}\int\limits_{0}^{|x|/||x||_L} t^{n-2} f_{n-1}\left(\frac{tx}{|x|}\right)dt\right)^\wedge(\xi).
$$
Integrate the latter inequality over $S^{n-1}$ and apply the spherical Parseval identity, Proposition \ref{procont}. Finally, the same computations (based on Lemma \ref{ell}) as in the proof of Theorem \ref{th:maina} give
$$
\mu_n(D)> \mu_n(L).
$$
\endpf

\section{Applications}\label{sec:apl}

\begin{co}\label{co:me}
Assume  $f_n(x)=f_{n-1}(x)$, then the answer to the BPGM problem  is 
affirmative if $n\le 4$ and negative if $n\ge 5.$ 
\end{co}
\pf In this case $tf_n(tx)/f_{n-1}(tx)=t$ is an increasing function, so we 
may apply  Theorems \ref{th:maina} and \ref{th:mainn}. First note that
$$
\|x\|_K^{-1}\frac{f_n\left(\frac{x}{\|x\|_K}\right)}
{f_{n-1}\left(\frac{x}{\|x\|_K}\right)}=
||\cdot||_K^{-1}.
$$ 
Thus  we may use the fact that for any origin symmetric convex body $K$ in 
$\R^n$, $n\le 4$,  
$||\cdot||_K^{-1}$ represents a positive definite distribution 
(see \cite{GKS}, \cite{K8}) to give the affirmative answer to BPGM in 
this case. 

For  $n\ge 5$, we first note that there is an infinitely smooth, symmetric, convex body
$L\subset \R^n$ with positive curvature and such that
 $||\cdot||_L^{-1}$ is not positive definite 
(see \cite{GKS}, \cite{K8}) and thus
$$
\|x\|_L^{-1}\frac{f_n\left(\frac{x}{\|x\|_L}\right)}
{f_{n-1}\left(\frac{x}{\|x\|_L}\right)}=
||\cdot||_L^{-1}
$$ 
is in $C^\infty(\R^n\setminus\{0\})$ and does not represent a positive definite distribution.

 Finally if $f_{n-1}\not\in C^\infty(\R^n\setminus\{0\})$ we finish the proof
 by approximating $f_{n-1}$ (and thus $f_n$)  by a sequence of strictly 
 positive functions belonging to $C^\infty(\R^n\setminus\{0\})$ . 
\endpf
\noindent {\bf Remark:} Note that the answers for the original Busemann-Petty problem and the Busemann-Petty problem for Gaussian Measures are particular cases of Corollary \ref{co:me},  with $f_n(x)=1$ and $f_n(x)=e^{-|x|^2/2}$ respectively.

\begin{lemma}\label{lm:1} Consider a symmetric star-shaped body $M \subset \R^n$ such that
 $\|x\|_M^{-1}$ is positive definite, then for any symmetric star-shaped bodies 
$K, L\subset \R^n$ such that, for every $\xi\in S^{n-1}$,
\begin{equation}
\int_{K\cap \xi^\bot}\|x\|_M dx \le \int_{L\cap \xi^\bot}\|x\|_M dx 
\end{equation}  
we have
$$  
\vol_n(K)\le \vol_n(L).  
$$  
\end{lemma}
\pf This follows from 
Theorem \ref{th:maina} with  $f_n=1$, $f_{n-1}=\|x\|_M$. In this case,
$$
\|x\|_K^{-1}\frac{f_n\left(\frac{x}{\|x\|_K}\right)}
{f_{n-1}\left(\frac{x}{\|x\|_K}\right)}=\|x\|_K^{-1}\frac{1}
{\|\frac{x}{\|x\|_K}\|_M}=||\cdot||_M^{-1}.
$$
\endpf

\begin{lemma}\label{lm:-1} Consider a symmetric star-shaped 
body $M \subset \R^n$ such that
 $\|x\|_M^{-1}$ is positive definite. Then for any symmetric star-shaped bodies 
$K, L\subset \R^n$ such that, for every $\xi\in S^{n-1}$,
\begin{equation}
\vol_{n-1}(K\cap \xi^\bot)\le \vol_{n-1}(L\cap \xi^\bot), 
\end{equation}  
we have 
$$
  \int_{K}\|x\|_M^{-1} dx \le \int_{L}\|x\|_M^{-1} dx 
$$  
\end{lemma}
\pf This theorem follows by the same argument as in Lemma \ref{lm:1}, but with the functions  $f_n=\|x\|_M^{-1}$, $f_{n-1}=1$.
\endpf
\noindent{\bf Remark: } It follows from  Theorem 
\ref{th:mainn}, and the standard approximation argument, that Lemmas \ref{lm:1} and \ref{lm:-1} are not true (even, with  additional convexity assumption) if  $\|x\|_M^{-1}$ is not positive definite.

Next we will use Theorem \ref{th:maina} and ideas from 
Lemmas \ref{lm:1} and \ref{lm:-1} to 
show some results on a lower dimensional version of the BPGM.
\begin{theorem}\label{th:codim} Consider a symmetric star-shaped body $M \subset \R^n$ such that
 $\|x\|_M^{-1}$ is positive definite, then for any symmetric star-shaped bodies 
$K, L\subset \R^n$, and $1\le k <n$ such that, for every $H\in G(n,n-k)$
$$
\int_{K\cap H}||x||_M^k dx \le \int_{K\cap H}||x||_M^k dx,
$$
we have
$$
\vol_n(K)\le \vol_n(L).
$$
\end{theorem}
\pf It was proved in \cite{GW} that every hyperplane section of intersection 
body 
is also an intersection body.  Using the 
relation of positive definite distributions to intersection
 bodies we get that if  $\|x\|_M^{-1}$ is a positive definite 
distribution than the restriction of  $\|x\|_M^{-1}$ on subspace $F$, is also
a positive definite distribution. Thus we may apply Theorem \ref{th:maina}
with functions $f_{i-1}(x)=||x||_M^{n-i+1}$ and $f_{i}(x)=||x||_M^{n-i}$.
Indeed, in this case
$$
\|x\|_K^{-1}\frac{f_i\left(\frac{x}{\|x\|_K}\right)}
{f_{i-1}\left(\frac{x}{\|x\|_K}\right)}=\|x\|_K^{-1}\frac{\|\frac{x}{\|x\|_K}\|_M^{n-i}}
{\|\frac{x}{\|x\|_K}\|_M^{n-i+1}}=||\cdot||_M^{-1}
$$
is a positive definite distribution. So, from Theorem \ref{th:maina}, we get that if for every $H\in G(n,n-i)$
$$
\int_{K\cap H}||x||_M^i dx \le \int_{K\cap H}||x||_M^i dx,
$$
then for every $F\in G(n,n-i+1)$
$$
\int_{K\cap F}||x||_M^{i-1} dx \le \int_{K\cap F}||x||_M^{i-1} dx.
$$
We iterate this procedure for $i=k,k-1,\dots,1$ to finish the proof.
\endpf
We can present a different version of  Theorem \ref{th:codim}, in a special cases
of $n-k=\{2,3\}$ and $K$ is a convex symmetric body:
\begin{co}\label{co:krz} Consider a symmetric star-shaped body $M \subset \R^n$, $n\ge 4$,
 such that
 $\|x\|_M^{-1}$ is positive definite. Fix $k$ such that $n-k\in\{2,3\}$, 
then for any convex symmetric  bodies $K, L\subset \R^n$, such that,
for every $H\in G(n,n-k)$
\begin{equation}\label{eq:krz}
\int_{K\cap H}||x||_M^{n-4} dx \le \int_{K\cap H}||x||_M^{n-4} dx,
\end{equation}
we have
$$
\vol_n(K)\le \vol_n(L).
$$
\end{co}
\pf We use the same iteration procedure as in Theorem \ref{th:codim}, 
but the first
steps of iteration procedure, for example,
 from subspaces of dimension $2$ to subspaces  of dimension $3$ is different and we  use Corollary \ref{co:me} with $f_2(x)=f_3(x)=||x||_M^{n-4}$. We use the
same idea to iterate from dimension $3$ to $4$.
\endpf
\noindent{\bf Remark:} Note that if $n-k=2$, then Corollary \ref{co:krz} is still
true with power $n-3$, instead of $n-4$ in (\ref{eq:krz}).
Corollary \ref{co:krz} is a generalization
of a result of Koldobsky (\cite{K6}, Theorem 8; see also \cite{RZ}), where
the case of $n-k=3$ and $M=B_2^n$ was considered.

 We also note that the generalization of the standard Busemann-Petty problem is open for those dimensions (see \cite{BZ}, \cite{K6}, \cite{RZ}).

Let $\|x\|_p=\left(\sum_{i=1}^n|x_i|^p\right)^{1/p}$, 
$\|x\|_\infty=\max\limits_{i}|x_i|$  and $B_p^n=\{x\in\R^n: \|x\|_p\le 1\}$.
\begin{co}\label{co:p} Consider $p\in (0,2]$ and $n\ge 2$. Then for any symmetric star-shaped bodies 
$K, L\subset \R^n$ such that
\begin{equation}
\int_{K\cap \xi^\bot}\|x\|_p dx \le \int_{L\cap \xi^\bot}\|x\|_p dx 
\end{equation}  
for every $\xi\in S^{n-1}$, we have
$$  
\vol_n(K)\le \vol_n(L).  
$$  
\end{co}
\pf This follows from Lemma  \ref{lm:1} and the fact that $\|x\|_p^{-1}$
is positive definite for $p \in (0,2]$ (see \cite{K3}, \cite{K8}).
\endpf
\noindent{\bf Remark: } Note that $\|x\|_p^{-1}$ does not represent a positive 
 definite distribution when
$p\in (2,\infty]$ and $n>4$ (see \cite{K3}, \cite{K8}), 
thus applying Theorem \ref{th:mainn}, together with the standard 
approximation argument, we get that
the statement of Corollary \ref{co:p} is not true in those cases,
 and counterexamples (even with bodies  $K$ and $L$ being convex) can be
 constructed.

Next we would like to use Lemma \ref{lm:1} to give a lower bound
on the integral over  a hyperplane section of the convex body.

\begin{lemma}\label{lm:int} 
Consider a symmetric star-shaped 
body $M \subset \R^n$ such that
 $\|x\|_M^{-1}$ is positive definite. Then for
any star-shaped body $K \subset \R^n$ there exits a direction
$\xi\in S^{n-1}$ such that
$$
\int_{K\cap \xi^\bot}\|x\|_M dx \ge \frac{n-1}{n}
\frac{\vol_{n-1}(M\cap\xi^\bot)}{\vol_n(M)}\vol_{n}(K).
$$
\end{lemma}
\pf Assume it is not true, then
$$
\int_{K\cap \xi^\bot}\|x\|_M dx <
\frac{n-1}{n}\frac{\vol_{n-1}(M\cap\xi^\bot)}{\vol_n(M)}\vol_{n}(K),\,\,\, \forall \xi\in
 S^{n-1}.
$$
Also note  that for $L\subset \R^{n-1}$
$$
\int_{L}\|x\|_L dx = \frac{1}{n}\int_{S^{n-1}} \|\theta\|_L^{-n}d\theta=\frac{n-1}{n}\vol_{n-1}(L).
$$ 
Applying the latter equality to $L=M\cap\xi^\bot$ to get
$$
\int_{K\cap \xi^\bot}\|x\|_M dx <
\frac{\vol_{n}(K)}{\vol_n(M)}\int_{M\cap\xi^\bot}\|x\|_M dx,\,\,\, \forall \xi\in S^{n-1}.
$$
Let $r^n=\frac{\vol_{n}(K)}{\vol_n(M)}$, then
$$
\int_{K\cap \xi^\bot}\|x\|_{rM} dx <
\int_{rM\cap\xi^\bot}\|x\|_{rM} dx,\,\,\, \forall \xi\in S^{n-1}.
$$
Thus
$$
\vol_n(K)< \vol(rM),
$$
or
$$
\vol_n(K)< \frac{\vol_{n}(K)}{\vol_n(M)}\vol_n(M)
$$
which gives a contradiction.
\endpf
\begin{co}\label{co:lp} For any $p\in [1,2]$ and any 
symmetric star-shaped body $K$ in $\R^n$ there exists
a direction $\xi \in S^{n-1}$ such that
$$
\int_{K\cap \xi^\bot}\|x\|_pdx \ge c_p n^{1/p}\vol_{n}(K),
$$
where $c_p$ is a constant depending on $p$ only.
\end{co}
\pf Again $\|x\|_p^{-1}$, $p\in [1,2]$ is a positive definite distribution 
(see \cite{K3}, \cite{K8}) and thus we may apply Lemma \ref{lm:int} (or Corollary \ref{co:p}) to get 
that there exist $\xi\in S^{n-1}$ such that
$$
\int_{K\cap \xi^\bot}\|x\|_p dx \ge \frac{n-1}{n}
\frac{\vol_{n-1}(B_p^n\cap\xi^\bot)}{\vol_n(B_p^n)}\vol_{n}(K).
$$
Next we use  that $B_p^n$, $p\in [1,2]$ is in the isotropic position, 
and the isotropic constant $L_{B_p^n} \le c$ (see \cite{Sc}), thus the ratio 
of volume of different hyperplane sections is bounded by two universal
 constants  (see \cite{Bo1}; \cite{MP}, Corollary 3.2):
$$
c\le \frac{\vol_{n-1}(B_p^n\cap\xi^\bot)}{\vol_{n-1}(B_p^n\cap\nu^\bot)}
\le C,\,\,\, \forall \xi,\nu\in S^{n-1},
$$
and
$$
c\vol_{n-1}(B_p^{n-1})\le \vol_{n-1}(B_p^n\cap\xi^\bot)
\le C\vol_{n-1}(B_p^{n-1}),\,\,\, \forall \xi\in S^{n-1}.
$$
Applying
$$
\vol_{n}(B_p^n)=\frac{[2\Gamma(1+\frac{1}{p})]^n}{\Gamma(1+\frac{n}{p})},
$$
we get
$$
\frac{\vol_{n-1}(B_p^{n-1})}{\vol_{n}(B_p^n)}=\frac{1}{2\Gamma(1+\frac{1}{p})}
\frac{\Gamma(1+\frac{n}{p})}{\Gamma(1+\frac{n-1}{p})} \ge c_p n^{1/p}.
$$
\endpf
\begin{co} For any convex symmetric body $K\in \R^n$, there are vectors
$\xi, \nu \in S^{n-1}$  such that $\xi\not=\pm \nu$ and
$$
\vol_{n-1}(K\cap \xi^\bot)\ge c\sqrt{\vol_{n-2}(K\cap \{\xi, \nu \}^\bot)\vol_n(K)},
$$
where $\{\xi,\nu \}^\bot$ is the subspace of codimension 2, orthogonal to $\xi$ and $\nu$.
\end{co}
\pf From Corollary \ref{co:lp} (with $p=1$) we get that there exists
a direction $\xi \in S^{n-1}$ such that
\begin{equation}\label{eq:k}
\int_{K\cap \xi^\bot}\sum\limits_{i=1}^n|x_i| dx \ge c n\vol_{n}(K).
\end{equation}
From continuity of volume measure, we may assume that 
$\xi \not\in \{\pm e_i\}_{i=1}^n$, where $\{e_i\}_{i=1}^n$ is the standard basis in $\R^n$.
From the inequality (\ref{eq:k}) we get that there exists $j\in \{1, \dots,n\}$ such that
\begin{equation}\label{eq:k2}
\int_{K\cap \xi^\bot}|x_j| dx \ge c \vol_{n}(K).
\end{equation}
Next we  apply the "inverse Holder" inequality (\cite{MP}, Corollary 2.7;
see also \cite{GrM}): for any symmetric convex body $L$ in $\R^{n-1}$ and unit vector $\theta\in S^{n-2}$
$$
\int_L |x\cdot \theta| dx \le c \frac{\vol_{n-1}^2(L)}
{\vol_{n-2}(L\cap\theta^\bot)},
$$
Using the latter inequality for $L=K\cap \xi^\bot$ and $\theta =e_j$ we get 
\begin{equation}\label{eq:lp}
\int_{K\cap \xi^\bot}  |x_j| dx \le c \frac{\vol_{n-1}^2(K\cap \xi^\bot) }{\vol_{n-1}(K\cap \{\xi, e_j\}^\bot)}.
\end{equation}
We compare inequalities (\ref{eq:k2}) and (\ref{eq:lp}) to finish the proof.
\endpf

\end{document}